\documentclass{article}

\usepackage{arxiv}

\usepackage{natbib}
\usepackage{graphicx}
\usepackage{epstopdf, epsfig}
\usepackage{amssymb}
\usepackage{amsmath}
\usepackage{hyperref}
\usepackage[nameinlink, noabbrev]{cleveref}    
\usepackage{xcolor}
\usepackage{subfig}
\usepackage{float}
\usepackage[normalem]{ulem}
\usepackage{soul}
\usepackage{colortbl,xcolor}

\hypersetup{
    colorlinks=true,
    linkcolor=blue,
    citecolor=blue,
    filecolor=magenta,      
    urlcolor=cyan,
    hypertexnames=true
    }

\usepackage{booktabs}       
\usepackage{nicefrac}       
\usepackage{microtype}      

\usepackage{ulem}

\title{Towards high-fidelity wind farm layout optimization\\ using polynomial chaos expansion and Kriging model}

\author{Yi-Xiao Shao$^1$,
        Zhen-Fan Wang$^1$,
        Shine Win Naung$^2$,
        Kai Zhang$^1$\thanks{Corresponding author: kai.zhang@sjtu.edu.cn},
        Yufeng Yao$^2$,
        Dai Zhou$^{1,3}$\\
        $^1$State Key Laboratory of Ocean Engineering, School of Ocean and Civil Engineering,\\ Shanghai Jiao Tong University, Shanghai, 200240, China \\
        $^2$School of Engineering, University of the West of England, Bristol, BS16 1QY, United Kingdom \\
        $^3$Shenzhen Research Institute of Shanghai Jiao Tong University, Shenzhen, 518063, China
        }

\begin{document}

\maketitle
\begin{abstract}

This paper presents a wind farm layout optimization framework that integrates polynomial chaos expansion, a Kriging model, and the expected improvement algorithm. The proposed framework addresses the computational challenges associated with high-fidelity wind farm simulations by significantly reducing the number of function evaluations required for accurate annual energy production predictions.  The polynomial chaos expansion-based prediction method achieves exceptional accuracy with reduced computational cost for over 96\%, significantly lowering the expense of training the ensuing surrogate model.  The Kriging model, combined with a genetic algorithm, is used for surrogate-based optimization, achieving comparable performance to direct optimization at a much-reduced computational cost.  The integration of the expected improvement algorithm enhances the global optimization capability of the framework, allowing it to escape local optima and achieve results that are either nearly identical to or even outperform those obtained through direct optimization.  The feasibility of the polynomial chaos expansion-Kriging framework is demonstrated through four case studies, including the optimization of wind farms with 8, 16, and 32 turbines using low-fidelity wake models, and a high-fidelity case using computational fluid dynamics simulations.  The results show that the proposed framework is highly effective in optimizing wind farm layouts, significantly reducing computational costs while maintaining or improving the accuracy of annual energy production predictions.

\end{abstract}

\section{Introduction}    \label{sec:Introduction}

In recent years, wind power has become a reliable alternative energy source to conventional fossil fuels. As the demand for wind energy increases, the size of the wind turbines grows, and the scale of wind farms expands rapidly as well. However, the aerodynamic efficiency of these large-scale wind farms is plagued by the wake effects within the turbine clusters, which leads to reduced wind speed and increased turbulence intensity for the downstream turbines. Studies show that wake effects typically result in power losses of 10\% to 20\% in some large offshore wind farms \citep{barthelmie2009modelling}.
 
Numerous studies have examined wind farm layout optimization (WFLO) as a promising approach to mitigating wake effects during the initial design phase of wind farms \citep{reddy2020wind}. By strategically positioning turbines to minimize wake interactions, WFLO aims to enhance the overall aerodynamic performance and energy output of the wind farm. This optimization process involves considering various factors such as turbine spacing, wind direction, and local topography to reduce the impact of wake losses on downstream turbines \citep{dong2021intelligent}. As the size of wind farms increases and offshore wind projects become more prevalent, the role of WFLO in achieving sustainable, high-efficiency wind energy generation becomes even more critical \citep{thomas2023comparison}.
A key component of WFLO is to select appropriate simulation tools for evaluating wake losses.
Two commonly used methods are low-fidelity analytical wake models and high-fidelity computational fluid dynamics (CFD).

The analytical wake models are widely used in WFLO due to high computational efficiency. 
These models use simplified functions based on momentum conservation or empirical relationships to describe the wake characteristics of turbines.
Among the various wake models, the Jensen model \citep{jensen1983note} and Gaussian model \citep{larsen1988simple} are typical analytical wake models that are frequently used.
A linear wake expansion is presumed by the Jensen model, with the velocity deficit depending on the distance, wind speed, and induction. 
The Gaussian wake model offers a wake shape that is more realistic by assuming a normal distribution for the velocity deficit.
The wake interactions in these analytical wake models are accounted for using the wake superposition methods \citep{machefaux2015engineering},
Common approaches include linear addition \citep{lissaman1979energy}, quadratic summation rule \citep{katic1987simple} and a parabolic type of approach \citep{larsen2009solitary}. 
Overall, the low-fidelity models are renowned for their computational efficiency and have been widely applied in large-scale wind farm layout optimizations. 
However, despite their widespread application, these models oversimplify the intricate physics of wake dynamics and wake interactions, which can lead to limited accuracy in predicting annual energy production (AEP).

CFD provides a higher-fidelity approach in simulating wind turbine aerodynamics. 
Within CFD, wind turbine modeling can be broadly classified into three categories: the blade-resolved model, the actuator line model (ALM), and the actuator disk model (ADM) \citep{mittal2016blade,zhang2023comparative}. 
By including boundary layer effects, turbulence, and dynamic stall, blade-resolved models provide an in-depth understanding of the flow surrounding individual blades  \citep{liu2017establishing}; however, they demand highly refined grids in the vicinity of the blades, leading to prohibitively high computational costs \citep{de2022blade}. 
To mitigate computational expense, the actuator line method (ALM) employs a discretization of rotor blades into lines of actuator points \citep{sorensen2002numerical}. 
These points project lift and drag forces, derived from pre-tabulated airfoil performance data. The actuator disk method (ADM) further simplifies this representation by modeling the rotor's swept area as a permeable disk \citep{mikkelsen2004actuator,svenning2010implementation,sanderse2011review}. 
On this disk, torque and thrust are implemented as source terms, inducing a reduction in flow velocity and pressure.
These actuator-based methods are widely used in high-fidelity simulations for wind farm aerodynamics studies \citep{stevens2017flow}, offering valuable insights into wake interference effects within turbine clusters \citep{shapiro2022turbulence,stevens2018comparison}.

Although the actuator-based CFD simulations offer a balance of prediction accuracy and computational efficiency in high-fidelity wind farm simulations, their direct application in WFLO remains computationally infeasible. 
The optimization algorithms for WFLO are broadly categorized as gradient-based methods which utilize the gradient information of the objective function with respect to the design variables to iteratively move towards an optimal solution, and population-based methods which operate on a population of potential solutions, iteratively improving the population using stochastic operators \citep{martins2021engineering} .
Regardless of the specific optimization method employed, WFLO typically requires evaluating numerous potential layouts over tens of thousands of iterations \citep{antonini2020optimal}.
This intensive computational requirement has constrained most high-fidelity CFD-based WFLO studies to small-scale wind farms to maintain manageable computational resources \citep{cruz2020wind}.

Another limitation of these high-fidelity WFLO studies is that they rely on a coarse-grained wind rose.
The wind rose describes the probabilistic distribution of wind direction and velocity at a site, and is crucial for accurately estimating a wind farm's AEP. 
Traditional WFLO evaluates energy production for each combination of wind speed and direction within the wind rose; these values are then weighted by their occurrence probabilities to determine the overall AEP. 
Thus, an overly simplistic wind rose may fail to capture the nuances of the wind distribution, potentially leading to suboptimal wind farm layouts. 
Conversely, increasing the wind rose's resolution introduces significant computational challenges, as exhaustively evaluating all wind speed and direction combinations can become computationally infeasible, especially when high-fidelity models are used.

To solve these problems, this paper proposes an integrated approach that leverages uncertainty quantification (UQ) and surrogate-based optimization (SBO) methods to significantly reduce the number of function calls in WFLO, thus enabling the use of high-fidelity CFD simulations in wind farm design.
To address the uncertain wind direction and speed encompassed in realistic wind roses, polynomial chaos expansion (PCE) is employed as the UQ tool to efficiently calculate the AEP for each wind farm layout during each iteration of the optimization process \citep{padron2019polynomial,dong2024wind}.
Furthermore, the SBO technique, trained on a limited set of high-fidelity simulations, approximates the relationship between turbine layout parameters and wind farm performance, allowing the optimization algorithm to navigate the design space much more efficiently \citep{bempedelis2023data,wang2024optimization}.
The integration of PCE and SBO aims to provide a robust and computationally efficient framework for high-fidelity WFLO, addressing the dual challenges of resolving a detailed wind rose and the iterative nature of the optimization process while minimizing computational costs.
This approach represents an important step towards high-fidelity simulations for efficient wind farm design optimization.

\section{Methods}    \label{sec:method}
In this section, the methods applied in this study are outlined in detail. 
First, the WFLO problem is formulated in Section \ref{sec:WFLO}. 
The low- and high-fidelity wind farm simulation techniques are introduced in Section \ref{sec:simulation}. 
Next, Section \ref{sec:UQ} discusses the application of polynomial chaos expansion (PCE) to calculate the AEP of the wind farm. 
The surrogate-based optimization methods are then presented in Section \ref{sec:SBO}, including the Kriging model, the expected improvement (EI) algorithm, and the genetic algorithm. 
Finally, the wind farm layout optimization framework is summarized in Section \ref{sec:framework}.

\subsection{Formulation of WFLO problem}
\label{sec:WFLO}
WFLO aims to maximize the AEP by optimizing wind turbine positions.
The wind farm's AEP is calculated by summing the power output of each turbine across various wind conditions (different wind speeds and directions), weighted by how frequently each condition occurs:
\begin{equation}
    \mathrm{AEP}=8760\sum_{i=1}^{n}\sum_{j=1}^{m}P_{j}(x_{i},y_{i})f_{j},
\end{equation}
where $P_{j}(x_{i},y_{i})$ is the power output of the $i$-th wind turbine located at $(x_i, y_i)$ under wind condition $j$, with $f_{j}$ representing its frequency of occurrence. 
The wind turbines are only allowed to be installed within the pre-specified wind farm boundary for the WFLO problem. 
The constraints for the WFLO problem are therefore expressed as follows:
\begin{equation}
    {\rm subject \ to}\  (x_i,y_i) \in \boldsymbol{D},
\end{equation}
where $\boldsymbol{D}$ denotes the boundary of the wind farm.
In addition, any two turbines are separated by a minimum distance of two diameters of the rotor to avoid blade collisions.

\subsection{Wind farm simulation approaches}
\label{sec:simulation}
\subsubsection{Low-fidelity model}
Although the aim of this study is to enable high-fidelity WFLO, the low-fidelity Gaussian wake model is first used to demonstrate the effectiveness of the proposed framework.
The Gaussian wake model \citep{bastankhah2014new} assumes that the velocity deficit in the wake of a wind turbine follows a normal distribution. 
The analytical expression for the wake velocity loss of a single turbine is given by:
\begin{equation}
\label{eq:wakemodel}
\displaystyle{\frac{\triangle U}{U_\infty}=\left(1-\sqrt{1-\frac{C_T}{8(k^*x/d_0+\epsilon)^2}}\right)\times \exp\left(-\frac{1}{2(k^{*}x/d_0+\epsilon)^2} \left(\left(\frac{z-z_h}{d_0}\right)^2+\left(\frac{y}{d_0}\right)^2\right)\right)},
\end{equation}
where $U_{\infty}$ represents the free-stream wind velocity, and $C_T$ denotes the thrust coefficient. $x$, $y$ and $z$ correspond to the streamwise, spanwise and vertical coordinates, respectively, with $z_h$ denoting the hub height. $d_0$ represents the diameter of the wind turbine, and $k^*$ is the growth rate. The value of $\epsilon$ is determined by equating the total mass flow deficit rate at $x=0$. 
Besides, for wind farm, a linear superposition method is employed to model the combined wake effects of multiple turbines.
The Gaussian model is integrated into the open-source package FLORIS (FLOw Redirection and Induction in Steadystate) \citep{floris2021flow}, which is a comprehensive and widely used wind farm simulation software that integrates multiple steady-state wake models to simulate the wake interactions between turbines.

\subsubsection{High-fidelity model}
For the high-fidelity wind farm simulation, the wake interactions are modeled employing CFD simulations, with the rotor of the turbine represented as an actuator disk. 
Flow fields around the wind turbines are numerically resolved by solving the Reynolds-Averaged Navier-Stokes (RANS) equations, which are formulated as
\begin{subequations}
\label{eq:RANS}
\begin{flalign}
    \frac{\partial \overline{u}_i}{\partial x_i}&=0,\\
    \rho\frac{\partial \overline{u}_i}{\partial t}+\rho\frac{\partial(\overline{u}_i\overline{u}_j)}{\partial x_j}&=-\frac{\partial \overline{p}}{\partial x_j}+\frac{\partial}{\partial x_j}(\mu(\frac{\partial \overline{u}_i}{\partial x_j}+\frac{\partial \overline{u}_j}{\partial x_i})-\rho\overline{u'_j u'_i})+f_i,
\end{flalign}
\end{subequations}
where $x_i$ represents the Cartesian space coordinate, and $\overline{p}$ and $\overline{u_i}$ denote the time-averaged pressure and velocity. The variables $\mu$ and $\rho$ correspond to the dynamic viscosity and air density, respectively, while $f_i$ symbolizes the source term that the actuator disk introduces. Arising from the time-averaging process, the Reynolds stress is given by $\overline{\rho u'_j u'_i}$, and is modeled using the $k-\epsilon$ turbulence model. This introduces two new variables: the turbulent kinetic energy $k$ and the dissipation rate $\epsilon$:
\begin{subequations}
\label{eq:turbulence_model1}
\begin{align}
    \rho\frac{\partial k}{\partial t}+\rho\frac{\partial(\overline{u}_i k)}{\partial x_i}&=\frac{\partial}{\partial x_j}(\frac{\mu_t}{\sigma_k}\frac{\partial k}{\partial x_j})+P_k-\rho\epsilon,\\
    \rho\frac{\partial \epsilon}{\partial t}+\rho\frac{\partial(\overline{u}_i \epsilon)}{\partial x_i}&=\frac{\partial}{\partial x_j}(\frac{\mu_t}{\sigma_\epsilon}\frac{\partial \epsilon}{\partial x_j})+C_{1 \epsilon}\frac{\epsilon}{k}P_k-C_{2 \epsilon}\frac{{\epsilon}^2}{k}\rho,
\end{align}
\end{subequations}
with
\begin{subequations}
\label{eq:turbulence_model2}
\begin{align}
      &P_k=-\rho\overline{u'_j u'_i}\frac{\partial u_j}{\partial x_i}, \\
      &\mu_t=-\rho C_{\mu}\frac{k^2}{\epsilon},
\end{align}
\end{subequations}
where $C_{\mu}$, $\sigma_k$, $\sigma_\epsilon$, $C_{1 \epsilon}$, and $C_{2 \epsilon}$ are the five constants in the $k-{\epsilon}$ turbulence model.

The volume force term $f_i$ in equation \ref{eq:RANS}(b) is decomposed into the streamwise part $f_{ix}$ and rotational part $f_{i\theta}$, and are dispersed radially using the Goldstein optimal distribution \citep{goldstein1929vortex}:
\begin{subequations}
\label{eq:forces}
\begin{align}
    &f_{ix} = A_x r^* \sqrt{1 - r^*},\\
    &f_{i\theta} = A_\theta \frac{r^* \sqrt{1 - r^*}}{r^* (1 - r_h^\prime) + r_h^\prime},
\end{align}
\end{subequations}
with
\begin{subequations}
\begin{align}
    &r^* = \frac{r^\prime - r_h^\prime}{1 - r_h^\prime}, \quad r^\prime = \frac{r}{R_P}, \quad r_h^\prime = \frac{R_H}{R_P},\\
    &A_x = \frac{105}{8} \frac{T}{\pi \Delta (3R_H + 4R_P)(R_P - R_H)},\\
    &A_\theta = \frac{105}{8} \frac{Q}{\pi \Delta R_P (3R_P + 4R_H)(R_P - R_H)},
\end{align}
\end{subequations}
where $f_{ix}$ represents the axial force, $f_i$ is the tangential force, $r$ denotes the distance between the point and the disk center, $R_P$ and $R_H$ are the external and internal radii of the disk, respectively. 
The thrust of the rotor is denoted by $T$, the torque of it is $Q$, and the disk's thickness is represented by $\Delta$. 
These parameters are used to model the forces exerted by the actuator disk on the surrounding flow, which are essential for determining the interaction between the turbine and the wind field. 

The thrust force and torque shown in equations \ref{eq:forces} are required as input to the actuator disk model.
We adopt the methodology proposed by \citet{richmond2019evaluation}, which integrates the wind turbine's torque and thrust lookup table directly into the flow solver. 
In order to calculate the theoretical inflow velocity in the lookup table, an iterative approach based on the one-dimensional actuator disk theory \citep{burton2011wind} is implemented in the simpleFoam source code (OpenFOAM package). More details of this approach are provided in our previous study \citep{wang2024optimization}.
The power of each wind turbine is calculated by summing each cell's power output within the actuator disk regions.

\subsection{Polynomial chaos expansion for AEP calculation}
\label{sec:UQ}
We follow \citet{padron2019polynomial} and \citet{dong2024wind} to apply polynomial chaos expansion for efficient prediction of wind farm AEP under complex wind rose conditions. 
The core idea is to create a polynomial approximation of the power output of a wind farm as a function of uncertain variables (wind speed, wind direction, etc.). 

Let the wind farm power $P(\boldsymbol{\xi})$ be a function that depends on the uncertain variables $\boldsymbol{\xi}=(\xi_1,\xi_2,...)$. In this paper, we consider two uncertain variables, the wind direction ($\xi_1$) and wind speed ($\xi_2$). 
The wind farm power $P(\boldsymbol{\xi})$ is approximated using a polynomial expansion:
\begin{equation}
\label{eq:PCE}
P(\boldsymbol{\xi}) \approx  \sum_{i=0}^{p-1} \alpha_i \boldsymbol{\varphi}_i(\boldsymbol{\xi}),
\end{equation}
where $p-1$ represents the order of the polynomial expansion, which is selected through cross-validation \citep{hastie2009elements}, and is limited to 10 in this paper \citep{padron2019polynomial}.
$\alpha_i$ and $\boldsymbol{\varphi}_i(\boldsymbol{\xi})$ are the polynomial coefficients and orthogonal polynomial basis, respectively.
The expected value of the power is expressed as 
\begin{equation}
\label{eq:expectedP}
\mu_P = \mathbb{E}[P(\boldsymbol{\xi})] = \int_\Omega P(\boldsymbol{\xi})\rho(\boldsymbol{\xi})\mathrm{d}\boldsymbol{\xi} \ \approx \int_\Omega \sum_{i=0}^{p}\alpha_i\boldsymbol{\varphi}_i(\boldsymbol{\xi})\rho(\boldsymbol{\xi})\mathrm{d}\boldsymbol{\xi},
\end{equation}
where $\rho(\boldsymbol{\xi})$ is the joint probability density function (PDF) of the wind direction and speed in the domain $\Omega$.
In PCE, the polynomial bases are designed to be orthogonal with each other with respect to the PDFs, i.e.,
\begin{equation}
    \label{eq:orthogonality}
    \langle\boldsymbol{\varphi}_i(\boldsymbol{\xi}), \boldsymbol{\varphi}_j(\boldsymbol{\xi})\rangle = \int_\Omega \boldsymbol{\varphi}_i(\boldsymbol{\xi}) \boldsymbol{\varphi}_j(\boldsymbol{\xi})\rho(\boldsymbol{\xi})\mathrm{d}\boldsymbol{\xi} = \delta_{ij},
\end{equation}
where $\delta_{ij}$ is the Kronecker delta, which is 1 if $i=j$ and 0 otherwise.
Note that the zero-th polynomial basis function is $\boldsymbol{\boldsymbol{\varphi}}_0(\boldsymbol{\xi})=1$ and is orthogonal to all other polynomial basis functions. The last term in Equation (\ref{eq:expectedP}) is reduced to 
\begin{equation}
\label{eq:inner_product1}
\int_{\Omega} \alpha_i \boldsymbol{\varphi}_i(\boldsymbol{\xi}) \rho(\boldsymbol{\xi}) \mathrm{d}\boldsymbol{\xi} = \alpha_i \langle \boldsymbol{\varphi}_i(\boldsymbol{\xi}), \boldsymbol{\varphi}_0(\boldsymbol{\xi}) \rangle = 0 \quad \text{for} \quad i > 0,
\end{equation}
and
\begin{equation}
\label{eq:inner_product2}
\int_{\Omega} \alpha_0 \boldsymbol{\varphi}_0(\boldsymbol{\xi}) \rho(\boldsymbol{\xi}) \mathrm{d}\boldsymbol{\xi} = \alpha_i \langle \boldsymbol{\varphi}_0(\boldsymbol{\xi}), \boldsymbol{\varphi}_0(\boldsymbol{\xi}) \rangle = \alpha_0 \langle 1, 1 \rangle = \alpha_0 \quad \text{for} \quad i = 0.
\end{equation}
Therefore, the expected value of power is the zero-th order coefficient, 
\begin{equation}
    \mu_P \approx \alpha_0,
\end{equation}
and the AEP is calculated by $\textrm{AEP} = 8760 \times \mu_P$.

In this paper, the probability distributions of the wind directions and wind speed are considered independent from each other, thus their joint distribution is given by $\rho(\boldsymbol{\xi})=\prod_{i=1}^2 \rho_i(\xi_i)$. 
The probabilistic distributions of the wind speed and wind direction are derived from the wind roses.
The multidimensional basis functions $\boldsymbol{\varphi}_i(\boldsymbol{\xi})$ are given by products of one-dimensional orthogonal polynomials:
\begin{equation}
\label{eq:orthogonal_basis}
\boldsymbol{\varphi}_i(\boldsymbol{\xi}) = \prod_{j=1}^n \varphi_{i_j} (\xi_j).
\end{equation}
For the distributions of wind speed and direction derived empirically from wind conditions, classical orthogonal polynomials do not directly apply. 
Thus, the data-driven approach proposed by \citet{oladyshkin2012data} is used to generate custom orthogonal polynomial basis $\varphi(\xi)$.

The polynomial coefficients $\alpha_i$ are sought through a regression method by constructing a linear system: 
\begin{equation}
\label{eq:linear_system}
\boldsymbol{\Phi\alpha} = \boldsymbol{P},
\end{equation}
where $\boldsymbol{\alpha}= \begin{bmatrix}\alpha_0  & \alpha_1 & \cdots & \alpha_{p-1}\end{bmatrix}^{\mathrm{T}}$ is the coefficients of the polynomial of order $p-1$.
$\boldsymbol{P}= \begin{bmatrix}P^0 & P^1 & \cdots & P^{m}\end{bmatrix}^{\mathrm{T}}$ is the power obtained by carrying out $m$ times of wind farm simulations using the samples generated by the Latin hypercube sampling method in the $\xi_1$-$\xi_2$ space. 
$\boldsymbol{\Phi}=\begin{bmatrix}\boldsymbol{\varphi}_0(\boldsymbol{\xi}^1) & \boldsymbol{\varphi}_1(\boldsymbol{\xi}^1) & \cdots & \boldsymbol{\varphi}_{p-1}(\boldsymbol{\xi}^1) \\
\boldsymbol{\varphi}_0(\boldsymbol{\xi}^2) & \boldsymbol{\varphi}_1(\boldsymbol{\xi}^2) & \cdots & \boldsymbol{\varphi}_{p-1}(\boldsymbol{\xi}^2) \\
\vdots & \ddots & \vdots \\
\boldsymbol{\varphi}_0(\boldsymbol{\xi}^m) & \boldsymbol{\varphi}_1(\boldsymbol{\xi}^m) & \cdots & \boldsymbol{\varphi}_{p-1}(\boldsymbol{\xi}^m)\end{bmatrix}$ stores the values of the polynomial basis functions evaluated at the sample points.
Finally, the least squares method is used to determine the coefficients by minimizing the sum of squared residuals, thereby ensuring the best possible fit between the model and the observed data:
\begin{equation}
\label{eq:least_squares}
\boldsymbol{\alpha} = \arg \min \| \boldsymbol{\Phi \alpha} - \boldsymbol{P} \|_2^2.
\end{equation}

\subsection{Surrogate-based Optimization}
\label{sec:SBO}
\subsubsection{Kriging Model}
The Kriging model \citep{krige1951statistical} is used as the surrogate model for surrogate-based optimization (SBO). 
The specific form of the Kriging model is given by:
\begin{equation}
\label{eq:kriging}
\hat{f}(x) = g(x)^T \beta + \varepsilon(x),
\end{equation}
where $\hat{f}$ denotes the Kriging approximation, $x$ represents the input variable, and $g(x)^T$ is the polynomial basis function vector, which provides global approximation, $\beta$ is the polynomial regression coefficient, and the Gaussian process error function, represented by $\epsilon(x)$, offers the local approximation and has zero uncertainty at the training sample locations. The function $g(x)$ can be a constant, first-order, or second-order polynomial, although its form has minimal impact on accuracy. In this paper, the open-source code DAKOTA \citep{adams2020dakota} is used to establish the Kriging model, with a second-order polynomial basis function vector selected. The error term $\epsilon(x)$ approximates the discrepancies between the surrogate model and the real model, following a normal distribution. Its covariance matrix is obtained from the correlation function matrix: 
\begin{equation}
\label{eq:correlation_function}
\text{Cov}(\varepsilon(x^{(i)}), \varepsilon(x^{(j)})) = \sigma^2 r(x^{(i)}, x^{(j)}),
\end{equation}
where $r(x^{(i)}, x^{(j)})$ represents the correlation function of sample parameters $x^{(i)}$ and $x^{(j)}$.
The Latin Hypercube Sampling (LHS) \citep{mckay2000comparison} method is used to generate the initial samples to build the Kriging model.

\subsubsection{Acquisition function}
Surrogate models constructed from initial sample data are typically insufficient for locating global optima. 
Consequently, during the optimization process, new samples are iteratively incorporated to refine the surrogate model's accuracy. 
Generally, two distinct strategies are employed for sample selection: targeting either regions exhibiting potentially optimal solutions (exploitation) or those characterized by high predictive uncertainty (exploration). 
In our earlier WFLO research \citep{wang2024optimization}, we exclusively employed exploitation by feeding each iteration's best layout (i.e., maximum surrogate prediction, MSP) back into the surrogate model. While computationally efficient, this approach frequently converges to suboptimal solutions, due to the numerous local optima present in high-dimensional WFLO problems.

The acquisition functions (AF) handle the critical balance between exploitation and exploration. They are used within Bayesian optimization to propose the next sampling point, balancing known promising areas with unexplored, uncertain regions.
The expected improvement  criterion is one of the most widely used AFs. \citep{movckus1975bayesian,jones1998efficient,zhan2020expected}. 
EI works by calculating how much improvement we can expect over the best solution found so far.
Let $f(x)$ be the objective function AEP and let $f(x^+)$ be the maximum observed AEP of $f$ at the sampled points $x_1, x_2, ..., x_n$. The EI value at a new point $x$ is given by:
\begin{equation}
\label{eq:EI}
EI(x) = \mathbb{E}[\max(0, f(x) - f(x^+))].
\end{equation}
The EI is expressed analytically using the predicted mean $\mu(x)$ and standard deviation $\sigma(x)$ of the Kriging model:
\begin{equation}
\label{eq:EI_kriging}
\displaystyle{EI(x) = (\mu(x) - \hat{f}(x))\Phi\left(\frac{\mu(x) - \hat{f}(x)}{\sigma(x)}\right) + \sigma(x) \phi\left(\frac{\mu(x) - \hat{f}(x)}{\sigma(x)}\right)},
\end{equation}
where $\Phi(\cdot)$ and $\phi(\cdot)$ are the cumulative distribution function and probability density function of the standard normal distribution, respectively. $\hat{f}(x)$ is the Kriging prediction.
When a region has a high predicted mean value $\mu(x)$, the first term on the right-hand side of equation (\ref{eq:EI_kriging}) becomes high, suggesting that this region is a good candidate for further exploitation.
Conversely, a larger standard deviation $\sigma(x)$ leads to a larger second term in the EI function, driving the algorithm to explore regions with higher uncertainty.

Conventionally, a single acquisition function is used in the surrogate-based optimization as the infilling criterion.
This paper applies a novel approach by first employing the MSP as the acquisition (the same as our previous paper \citep{wang2024optimization}), aiming at fast identification of a set of local optima. 
This approach is then augmented with the EI module as a second acquisition function to search for regions with uncertainty, thus enhancing the global optimization capability.
Notably, the convergence criterion for the MSP module is that the intermediate optimal layout must already be represented in the sample set. 
The stopping condition of the EI module is that the maximum EI value obtained through the optimization iteration is lower than the specified threshold (set at 0.1 in this paper).

\subsubsection{Genetic Algorithm}
This study uses the genetic algorithm (GA) as the optimization algorithm for SBO. The GA simulates the natural evolutionary process through computer modeling, which includes population reproduction, selection of adaptive individuals, and the adoption of natural processes like mutation and crossover. In the optimization problem, design parameters are treated as genes and encoded into bit arrays or strings to represent chromosomes. The optimization process is driven by the fitness of the objective function, with a new population generated through selection, crossover, and mutation to replace the old population until convergence criteria are met. The basic steps of the GA are as follows:
\begin{enumerate}
\item Initial population generation: Randomly generate the initial population and compute the fitness function for each individual. The initial population size should be set carefully.
\item Selection: Eliminate individuals that do not meet the constraints or penalize them with a penalty coefficient, then select parent individuals for mating based on fitness.
\item Crossover and mutation: Perform crossover and mutation operations based on predefined rates to generate offspring and avoid local optima.
\item Check convergence: If the convergence criteria are met or the maximum number of calculations is exceeded, the process ends; otherwise, repeat Steps 2-4.
\end{enumerate}

\subsection{The PCE-aided SBO framework}
\label{sec:framework}
\begin{figure}
\centering
\includegraphics[width=0.8\textwidth]{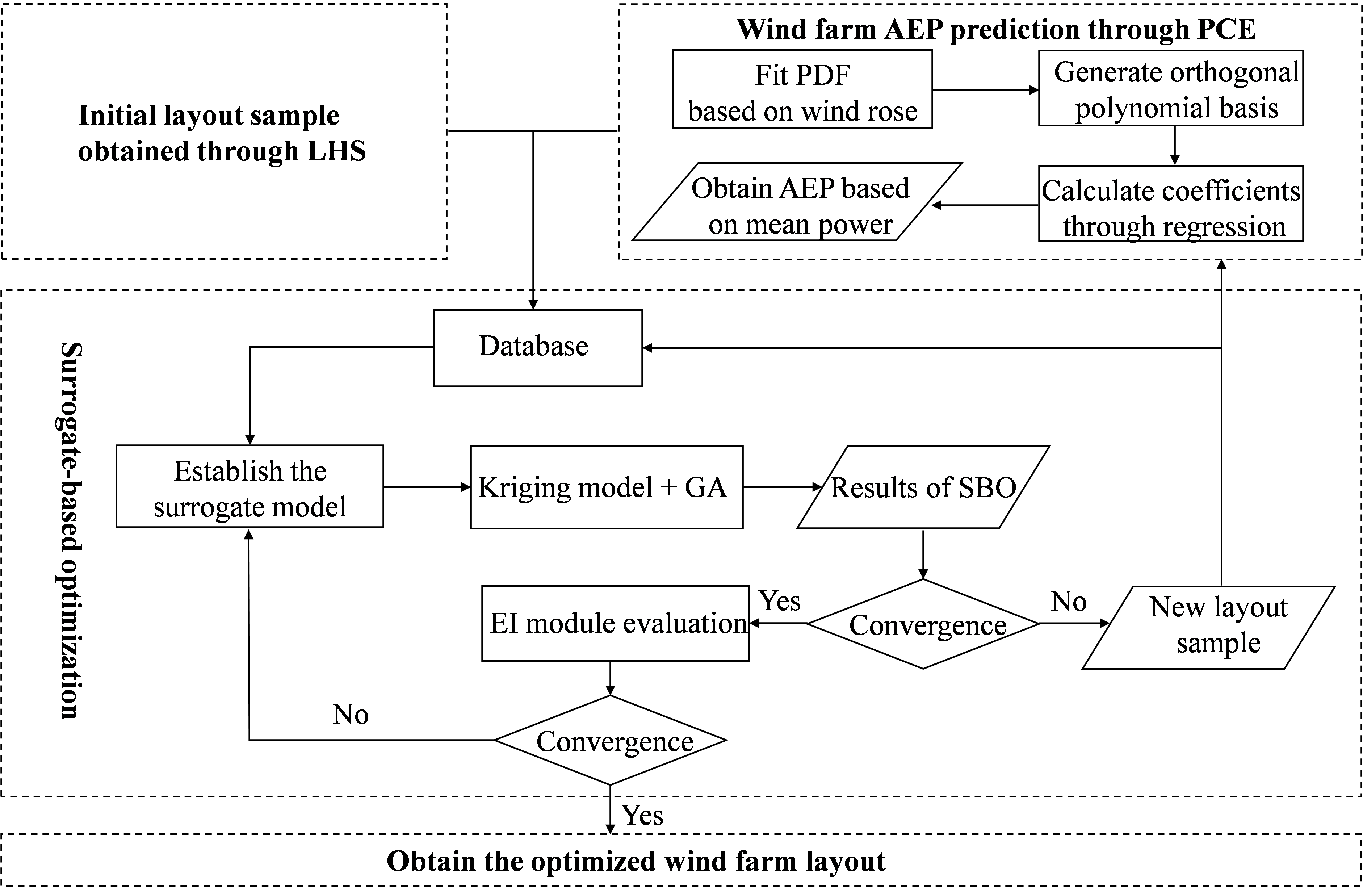}
\caption{The workflow of the adaptive optimization framework.}
\label{fig:workflow}
\end{figure}

Based on the methods outlined above, an adaptive optimization framework is developed in this study. The global surrogate model is integrated with the optimization algorithms. Additionally, the latest iteration, along with its true response, is added to the training dataset to enhance the accuracy of model predictions \citep{liu2021genetic}. The main steps of the framework are as follows: 
\begin{enumerate}
    \item Generate initial random layout dataset using Latin Hypercube Sampling (LHS);
    \item Apply the PCE model to calculate the AEP of the initial layout dataset;
    \item Construct the surrogate model (the Kriging model) based on the initial dataset;
    \item Integrate the surrogate model to the optimization algorithm to obtain the optimal results;
    \item Check if the optimal design has converged (the intermediate optimal layout must already be represented in the sample set) or if the maximum number of iterations has been exceeded. If not, repeat Steps 2–4. If converged, proceed to optimize the EI module;
    \item Check if the EI value satisfies the convergence condition (the maximum EI value obtained through the optimization iteration lower than the specified threshold, which is set to be 0.1). If not, repeat 2-5.
\end{enumerate}
The detailed flowchart of the framework is shown in Figure \ref{fig:workflow}. 

\section{Case Studies}
\label{sec:case}

\begin{figure}
\centering
\includegraphics[width=0.9\textwidth]{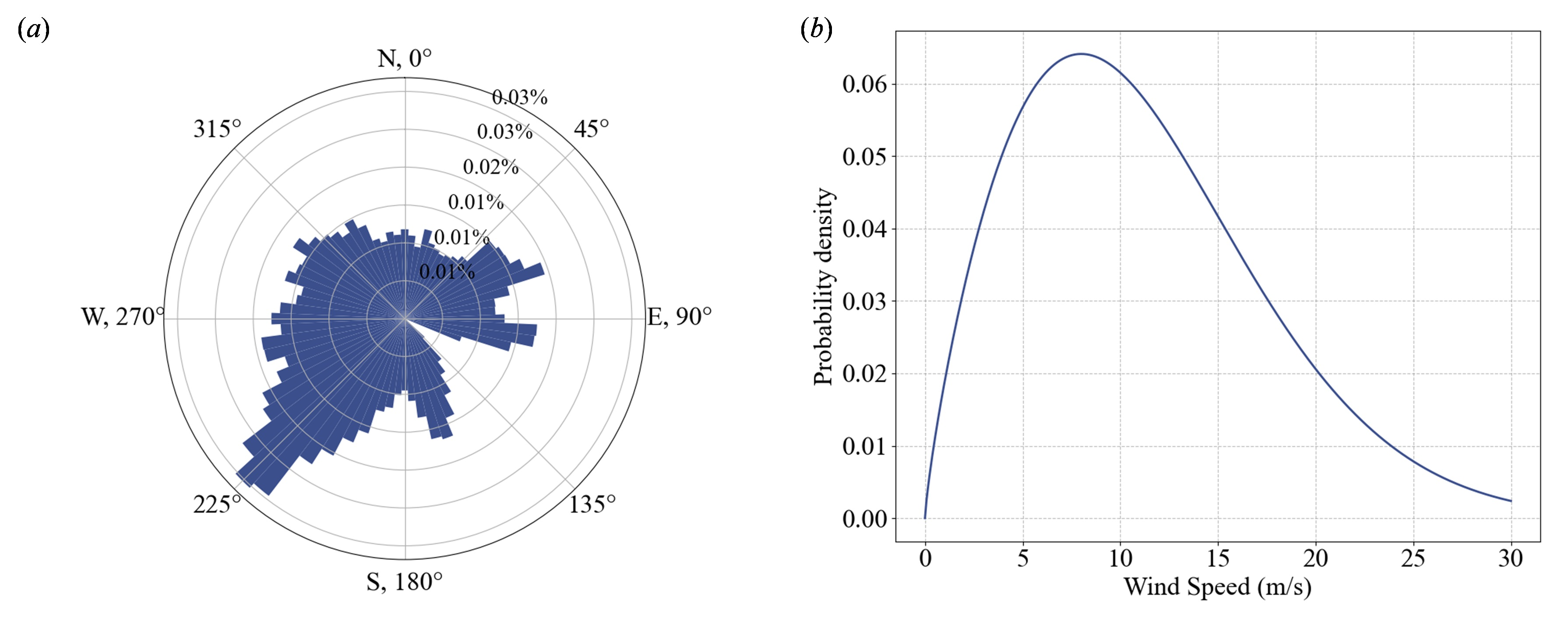}
\caption{($a$) The wind direction distribution with the bin width of 5 degree and ($b$) the Weibull distribution for wind speed. The distribution is truncated according to the cut-in speed (3 m/s) and cut-out speed (25 m/s). The bin width of wind speed is 1 m/s.}
\label{fig:windrose}
\end{figure}

In this section, we carry out four WFLO cases as shown in Table \ref{tab:case} to demonstrate the accuracy and effectiveness of the proposed framework.
In the first three cases, computationally efficient low-fidelity Gaussian wake models are employed to study the effects of a set of parameters on the optimization performance of the SBO approach.
For these cases, the wind data (figure \ref{fig:windrose}) sourced from \citet{gebraad2017maximization} is used, where the wind speed follows a Weibull distribution.
Next, the high-fidelity WFLO is demonstrated in the last case using CFD method.
In all cases, the NREL 5 MW reference turbine with a diameter $D=126$ m \citep{jonkman2009definition} is used as the wind turbine model.

\begin{table}
    \centering
    \caption{Setup for the four WFLO case studies. $D$ is the diameter of the NREL 5 MW wind turbine.}
    \begin{tabular}{cccccc}
        \toprule
        case & No. of turbines & fidelity &  wind farm region  & wind rose\\
        \midrule
        I & 8 & low-fidelity &  $8D\times8D$  & Figure \ref{fig:windrose} \\
        II & 16 & low-fidelity &  $12D\times12D$  & Figure \ref{fig:windrose}\\
        III & 32 & low-fidelity &  $16D\times16D$  & Figure \ref{fig:windrose}\\
        IV & 8 & high-fidelity &  $8D\times8D$  & Figure \ref{fig:windrose} with constant speed 8 m/s \\
        \bottomrule
    \end{tabular}
\label{tab:case}
\end{table}

\subsection{Effects of the sampling sizes}
For both PCE and Kriging, initial samples are required to train the respective models.
We evaluate the effects of sampling sizes on the performance of the  focus on case I, in which the positions of 8 wind turbines are to be optimized on a $8D\times 8D$ rectangular region.
Similar to our previous study \citep{wang2024optimization}, the wind turbines are only allowed to be located on the vertices of the $9\times 9$ evenly spaced grid, rendering it a discrete optimization problem.

\begin{figure}
\centering
\includegraphics[width=1.0\textwidth]{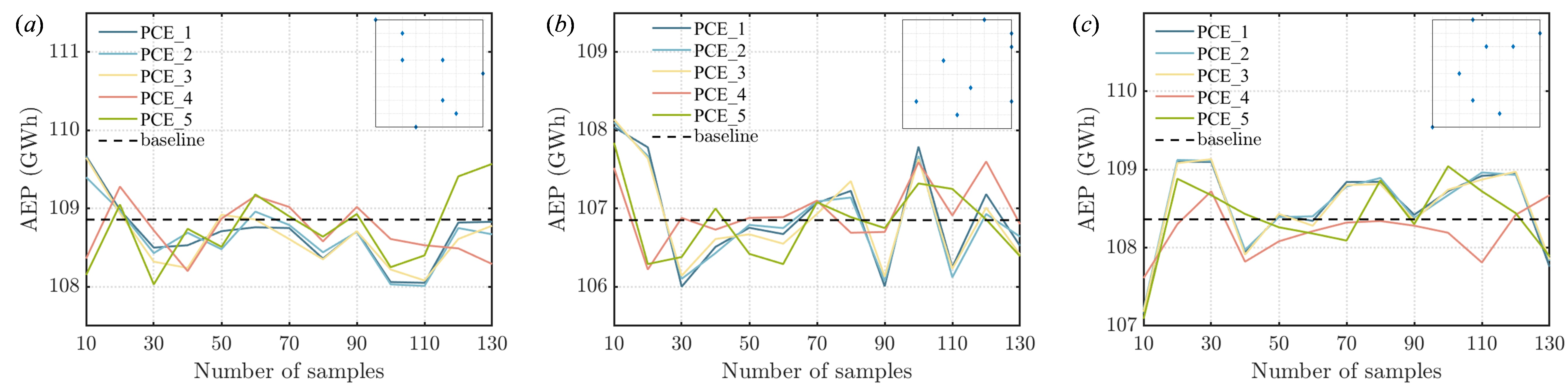}
\caption{Effects of sample size on the accuracy of AEP prediction in PCE. Three random layouts are tested, each with five instances of samples. The respective layouts are plotted in the upper right corner.}
\label{fig:error}
\end{figure}

The PCE module evaluates the AEP of an arbitrary wind farm layout by computing the polynomial coefficients $\boldsymbol{\alpha}$ through regression sampling of wind direction and speed, as described in equation \ref{eq:least_squares}.
Figure \ref{fig:error} illustrates the impact of the wind condition sample size on AEP accuracy using three random layouts. For each sample size, the AEP is computed 5 times using different random samples generated via LHS. 
For comparison, a baseline AEP is calculated by fully traversing the wind rose, which requires 1584 function evaluations (72 wind direction bins $\times$ 22 wind speed bins).
As shown in Figure \ref{fig:error}, the PCE-based AEP prediction method exhibits significant robustness, as evidenced by the consistent and reliable performance across multiple runs and layouts.
With increasing sample size, the AEPs of the three layouts predicted by PCE exhibit fluctuation around the baseline with an error of around 1\%, regardless of the instances of the wind condition samples.
For this study, a sample size of 50 is selected for all subsequent calculations, accounting for only 3.1\% of the computational demand associated with the conventional method.

\begin{figure}
\centering
\includegraphics[width=0.8\textwidth]{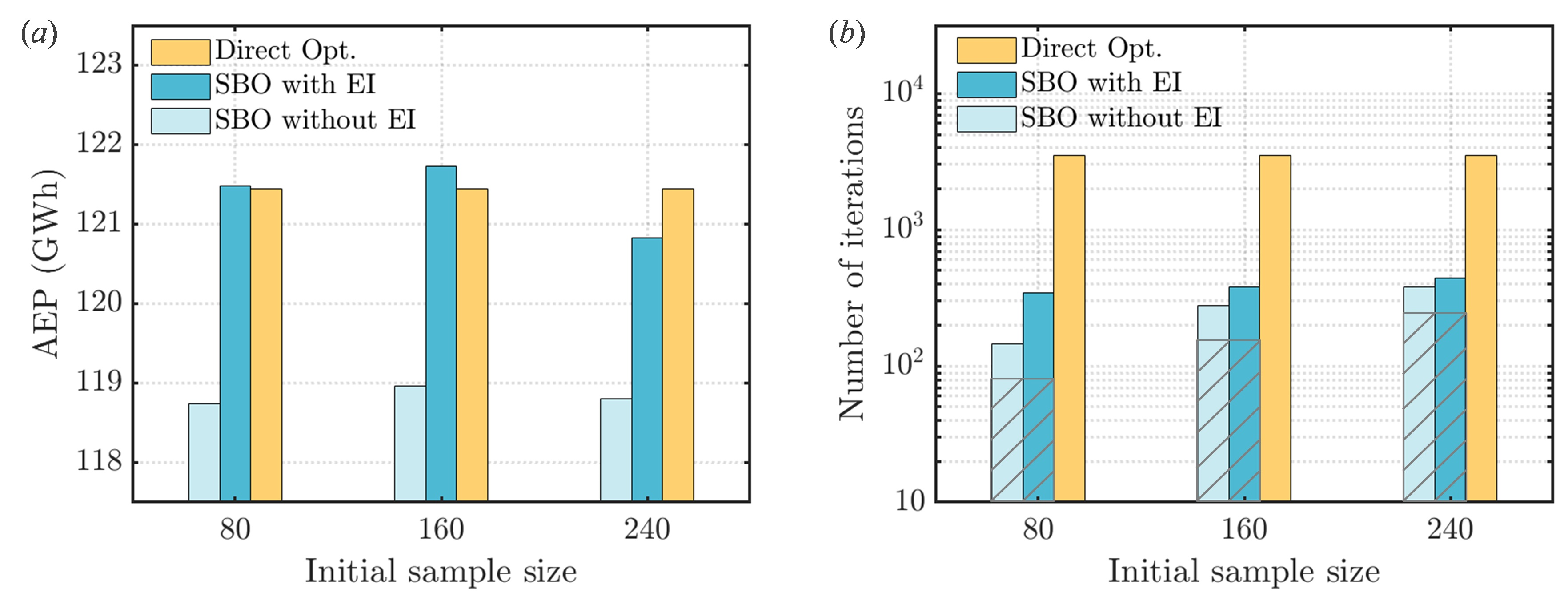}
\caption{Effects of the initial sample size on ($a$) the optimized AEP and  ($b$) number of iterations in the Kriging model. Note that in $(b)$, the number of iterations include the evaluations of initial sample size. The shaded area delineates the initial sample size.}
\label{fig:sample}
\end{figure}

Subsequently, the initial sample size of the wind farm layout designs for building the Kriging model is examined in Figure \ref{fig:sample}.
The initial samples are generated through LHS, and their sizes are set as 80, 160, and 240, corresponding to 5$d$, 10$d$, and 15$d$, where $d=16$ is the dimension of the optimization problem (8 wind turbines $\times$ 2 coordinates).
With growing initial sample size, the optimized AEPs through SBO with and without the EI module exhibit negligible variations.
While the total number of iterations (including initial training samples) required to reach the optimal solution increases with the initial sample size, this is primarily due to the larger initial sample set. 
The number of iterations during the optimization process itself actually decreases with increasing initial sample size. This suggests that a larger initial sample, potentially generated using parallel computing, can shorten the more computationally expensive optimization phase.
Compared to the direct optimization (traversing the wind rose) with a genetic algorithm as optimizer that requires $\sim$ 3000 iterations, the computational cost is significantly reduced by an order using SBO.
In case II and III, an initial sample size of $5d$ is used to carry out WFLO.

\begin{figure}
\centering
\includegraphics[width=0.7\textwidth]{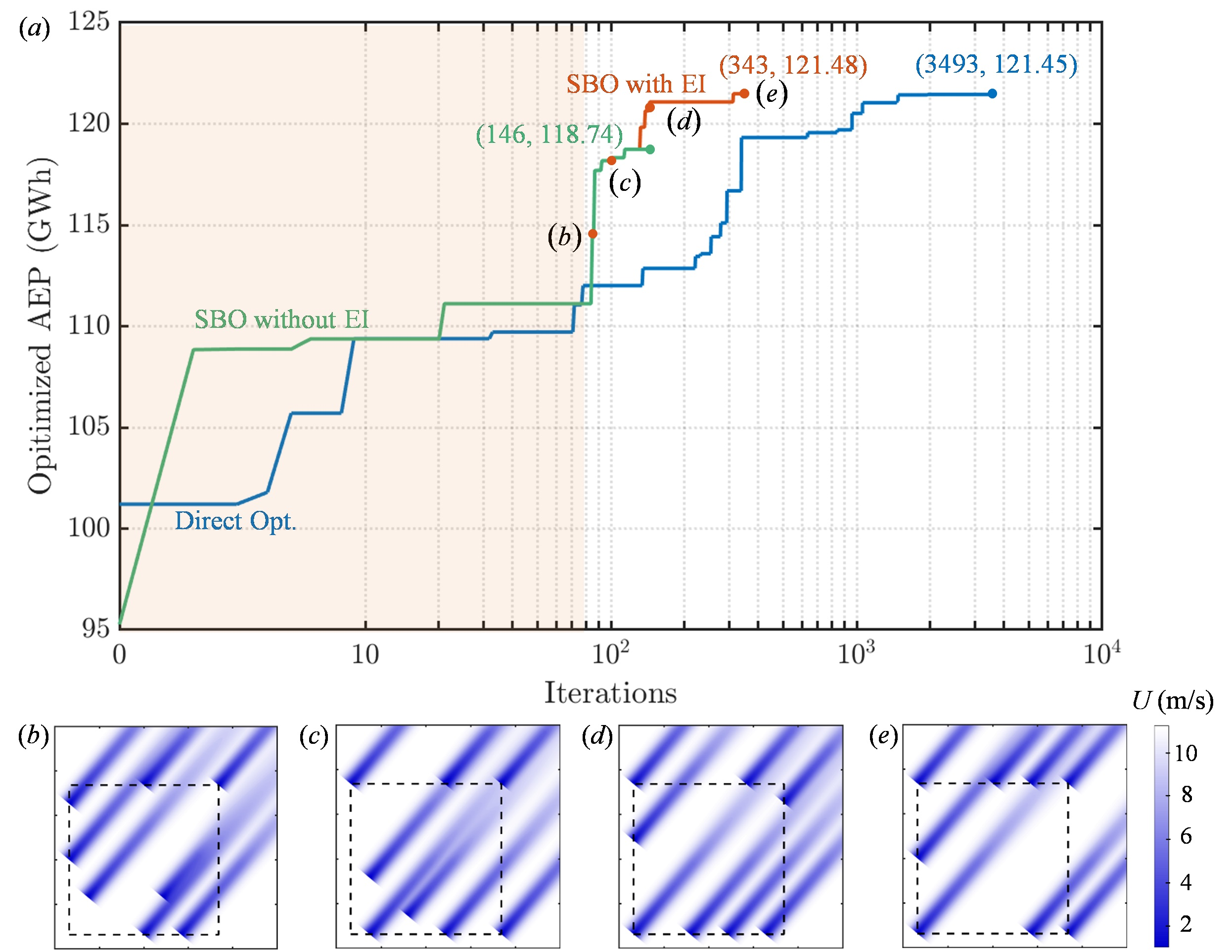}
\caption{($a$) Evolution of the optimized AEP using SBO (with and without EI) and direct optimization for case I, with intermediate optimal layout designs during the SBO (with EI) shown in $(b)-(d)$, and the final design in $(e)$. The shaded area in $(a)$ indicates the initial samples. The optimized AEPs for SBO with and without EI module follow the same path in the initial stage of the optimization. The dashed box in $(b)-(e)$ represents the boundary of wind farms.}
\label{fig:Iterations1}
\end{figure}

The impact of the EI module on optimization performance is also analyzed in Figure \ref{fig:sample}. 
Without the EI module, the SBO selects the point with the maximum predicted value from the Kriging model as the next sample point. 
Since the WFLO problem has the multi-modal nature \citep{wang2024optimization,thomas2023comparison}, this can lead to premature convergence to a local optimum. 
The EI module addresses this issue by guiding the search towards regions of higher uncertainty, allowing the optimization process to escape local optima.
This process is visualized in Figure \ref{fig:Iterations1}. 
Due to the architecture of the SBO framework as described in Figure \ref{fig:workflow}, both optimizations (with and without EI) follow the same trajectory in the initial iterations.
However, once the optimization without EI converges, the EI module takes over the optimization process to explore uncertain regions, continuing to improve the AEP. 
While this exploration requires additional computational resources, it clearly demonstrates the EI module's ability to enhance the global search capability of the framework and mitigate the risk of premature convergence.

\subsection{Application to larger wind farm cases}

In cases II and III, we use the proposed framework to determine the optimal arrangement of 16 and 32 wind turbines in $12D \times 12D$ and $16D \times 16D$ rectangular domains, respectively. 
The layout areas are discretized as $13 \times 13$ and $17 \times 17$ grids, totaling 169 and 289 vertices as potential turbine positions, respectively. 
Following case I, the PCE module uses 50 samples to compute the AEP of arbitrary wind farm layout design, and the initial sample size for constructing the Kriging model is fixed at $5d$, translating to 160 and 320 samples for cases II and III, respectively.
For both cases, the EI module is incorporated into the SBO.

\begin{figure}
\centering
\includegraphics[width=0.7\textwidth]{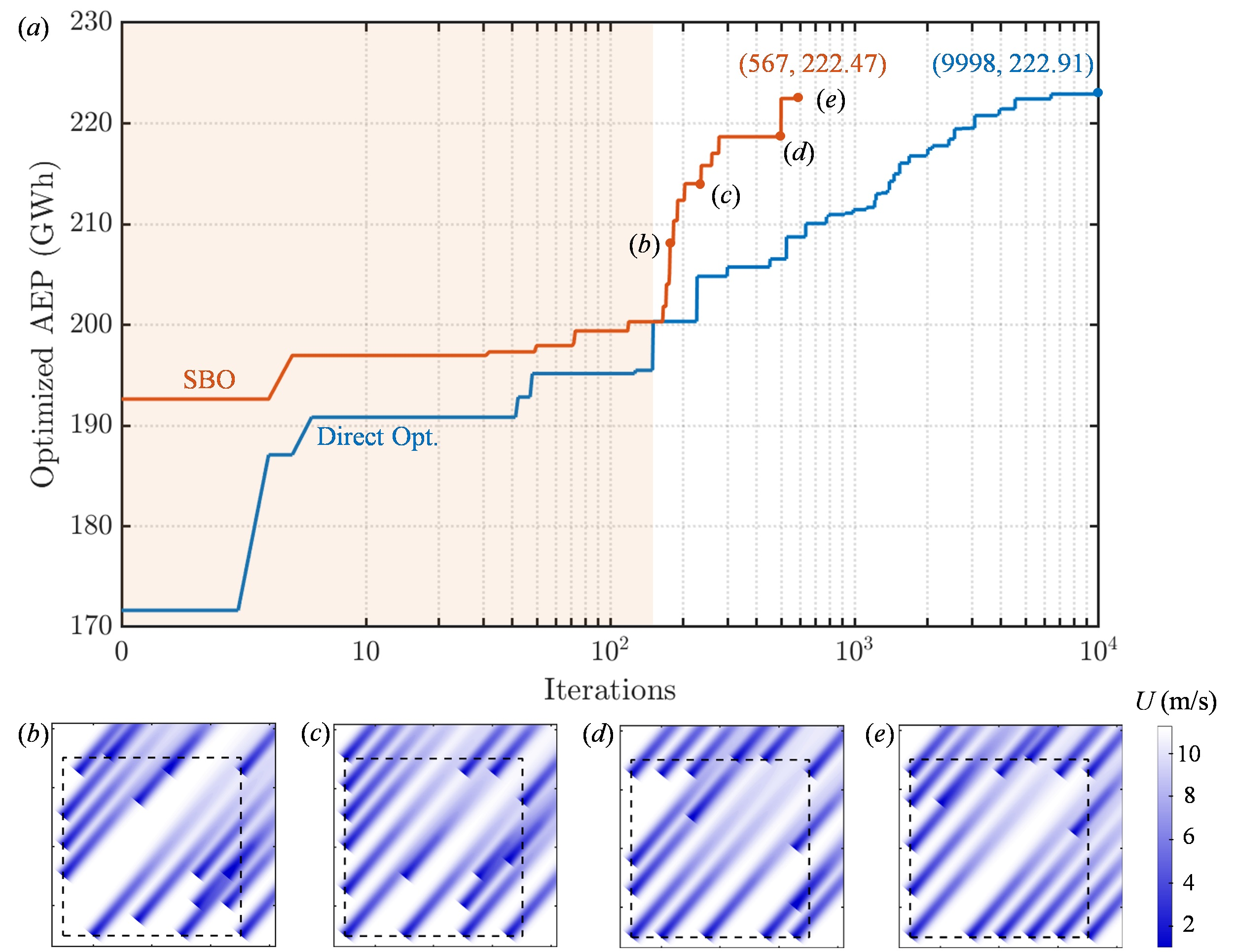}
\caption{Evolution of the optimized AEP using SBO (with EI module) and direct optimization for case II.}
\label{fig:Iterations2}
\end{figure}

\begin{figure}
\centering
\includegraphics[width=0.7\textwidth]{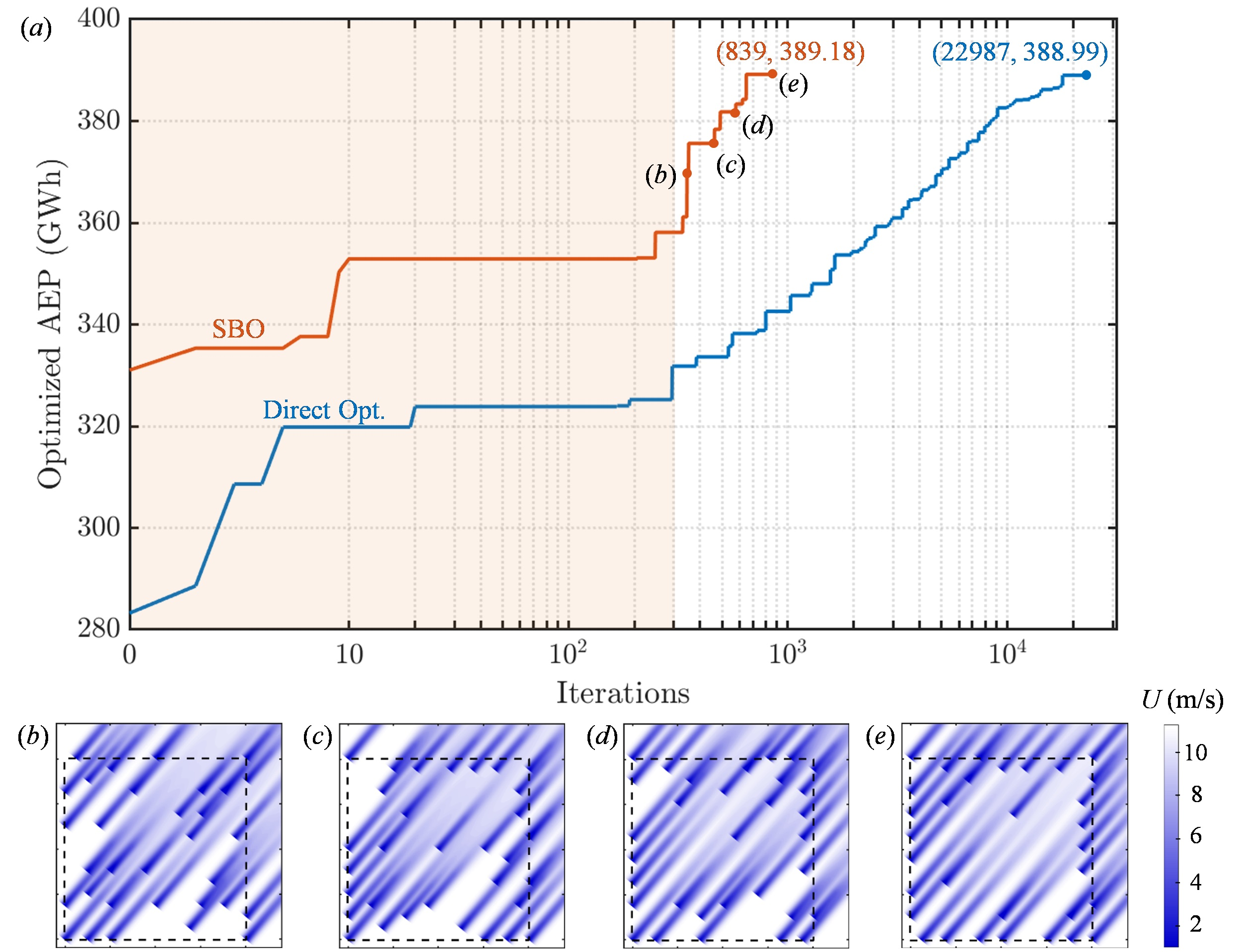}
\caption{Evolution of the optimized AEP using SBO (with EI module) and direct optimization for case III.}
\label{fig:Iterations3}
\end{figure}

\begin{figure}
    \centering\includegraphics[width=0.45\linewidth]{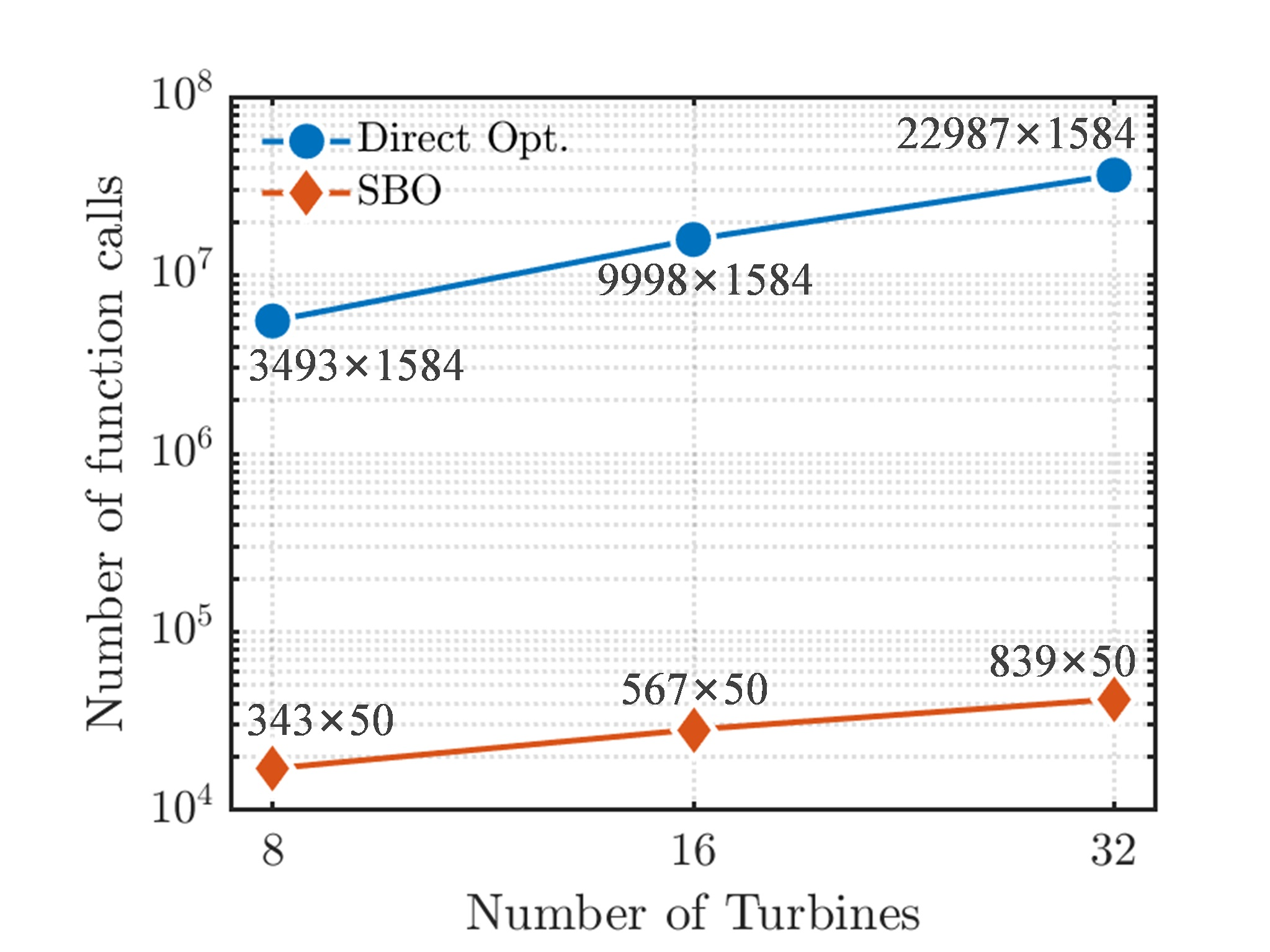}
    \caption{Comparison of function evaluations for SBO and direct optimization of WFLO. For direct optimization, the number of function calls is the product of the number of iterations in genetic algorithm until convergence and the number of wind speed and direction combinations. For SBO-PCE, it is calculated as the product of the total SBO iterations (including the initial sample) and the PCE sample size.}
    \label{fig:calc}
\end{figure}

Figures \ref{fig:Iterations2} and \ref{fig:Iterations3} depict the progression of the optimized AEP throughout the optimization process for cases II and III, respectively.
As the optimization advances, the turbines gradually migrate from the central regions of the wind farm domain to the periphery of the farm's boundary.
In both cases, the SBO approach converges to the optimal layout with substantially fewer iterations compared to direct optimization, achieving an AEP that is either very close to or surpasses that of the direct method.
As summarized in Figure \ref{fig:calc}, SBO requires only approximately 1/10, 1/20, and 3/100 of the iterations needed by direct optimization for the respective cases.
Furthermore, the PCE also provides substantial savings in the number of function calls during both the initial training phase and each iteration cycle.
Consequently, the PCE-aided SBO achieves comparable AEP as direct optimization with only $\sim$ 0.3\% of the total function calls.
This dramatic reduction in computational cost makes the application of high-fidelity models to WFLO problems far more feasible.

\subsection{Application to high-fidelity WFLO case}
\label{sec:ADM}
\begin{figure}
\centering
\includegraphics[width=0.8\textwidth]{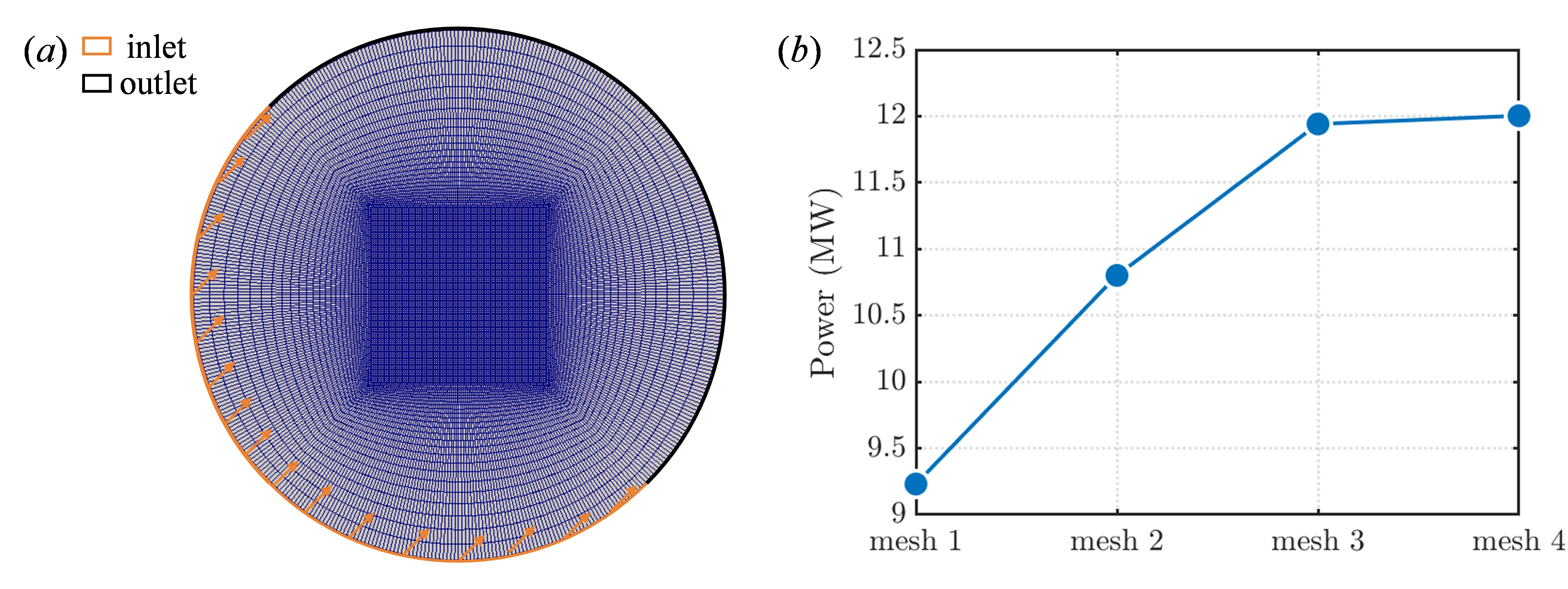}
\caption{($a$) The computational domain used in the CFD simulation and ($b$) the wind farm power obtained from different mesh designs.}
\label{fig:setup}
\end{figure}

In this section, we demonstrate the application of the proposed WFLO framework using high-fidelity CFD simulations in case IV.
Due to the limitation of our computational resources, case IV adopts the settings of case I to maximize the power output of a wind farm with 8 turbines in a rectangular space measuring $8D\times 8D$, with a further simplification of the wind speed to a constant value of 8 m/s (according to the weighted average of wind speed) for all wind directions.
For the same reason, the direct optimization as shown in the previous cases is not performed for this high-fidelity case.
As shown in Figure \ref{fig:setup}($a$), a circular computational domain is used to conveniently accommodate the wind velocity from different directions specified by the wind rose.
The radius of the circular computational domain is set as $12D$.
In the vertical direction, the height of the computational domain is set as 8$D$.
Besides, the atmospheric boundary layer at the inlet is modeled as $U_{in} = U_{\infty}(z/z_0)^{\alpha}$, where the reference velocity is set to $U_\infty=8$ m/s, the turbine hub height is $z_0 = 90 $ m, and the velocity profile's shear rate is fixed to $\alpha = 0.14$.
The outlet pressure is set to a reference pressure $p_\infty = 0$ Pa, with a zero-gradient condition for velocity.
A more detailed description of the case setup is provided in our previous study \citet{wang2024optimization}.

The native blockMesh tool in OpenFOAM is employed to create the mesh for the computational domain. The mesh is clustered within the wind farm boundary, with further refinement near the actuator disk zone to ensure accurate loading of the source terms for the ADM into the flow. 
A test of the mesh dependency is performed on a random layout to determine the optimal grid size. 
The minimum element sizes for meshes 1 through 4 are set to 0.1$D$, 0.05$D$, 0.03$D$ and 0.015$D$, respectively, resulting in a total of $1.06 \times 10^5$, $2.93 \times 10^5$, $6.26 \times 10^5$ and $9.72 \times 10^5$ control volumes. 
As shown in Figure \ref{fig:setup}($b$), meshes 1 and 2 significantly underestimate the power output. 
However, the difference between meshes 3 and 4 is negligible. 
Therefore, mesh 3 is selected for this study.
It is noted that one CFD simulation employing mesh 3 requires around 5 minutes to converge to the steady-state solution, using 36 processors with AMD EPYC 7302 CPU. 
In comparison, the execution of one evaluation of wind farm power using FLORIS barely takes more than 1 second. 

Since in case IV we only consider the wind direction as the uncertain variable, it is verified that only 8 wind direction samples are enough to obtain an accurate estimation of AEP.
In the discussion for case I, it is observed that increasing the initial sample size can shorten the number of iterations required in the optimization process.
Considering the fact that the generation of the initial layout samples can be parallelized but the optimization process cannot, we adopt a larger initial sample size of $10d=160$ for the high-dimensional WFLO.
Based on the Kriging model built on the 160 initial samples, the SBO eventually spends 112 iterations to find the optimal layout design.
All together, a total of 2176 ADM-based CFD simulations have been performed to find the optimal solution of the 8-turbine wind farm, under the wind rose specified in Figure \ref{fig:windrose} with a constant wind speed.

\begin{figure}
\centering
\includegraphics[width=0.8\textwidth]{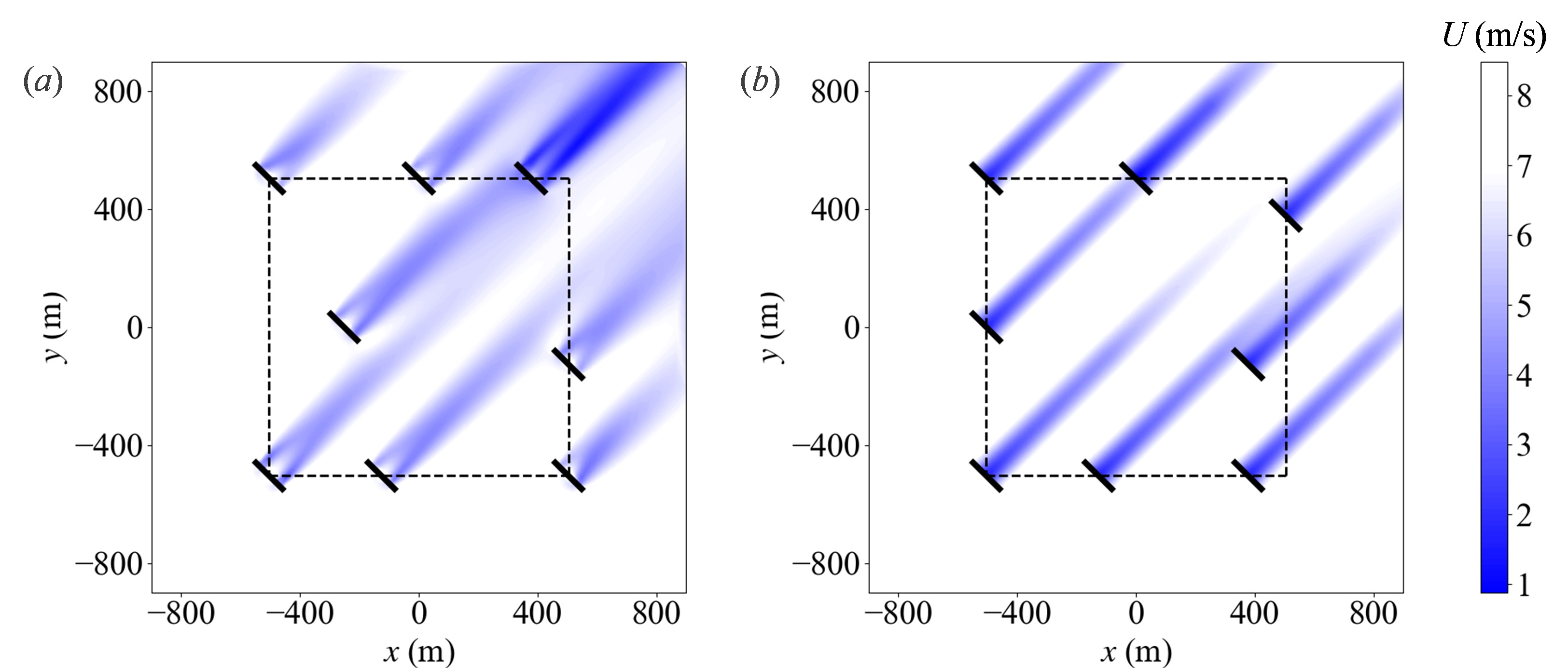}
\caption{The velocity field at hub height under the dominant wind direction of the optimal layout obtained from ($a$) CFD-based WFLO and ($b$) low-fidelity WFLO. The dashed box is the boundary of the wind farm.}
\label{fig:ADM_result}
\end{figure}

The velocity field of the optimal layout under the dominant wind direction is illustrated in Figure \ref{fig:ADM_result}($a$).
For comparison, the optimal layout obtained using the low-fidelity wake model, employing the same wind rose and SBO settings, is also presented in Figure \ref{fig:ADM_result}($b$).
It is clear that the CFD simulations are able to reproduce more intricate wake features such as the higher velocity region in the root region of the rotor, and wake asymmetry due to the mutual interactions among the turbines.
Despite the different layouts obtained from the two fidelities, the turbines are primarily distributed along the boundary in both cases. 
Furthermore, most of the turbines are delicately positioned in a staggered manner, or in the far wake region of the upstream turbine, to avoid wake interference effects.
Evaluated by the high-fidelity CFD simulation, the optimal layout obtained from the low-fidelity model is 106.79 MW, which is lower than the layout obtained by the CFD-based WFLO (108.51 MW).
The successful application of the PCE-Kriging framework in this high-fidelity CFD case demonstrates its potential to enable high-fidelity wind farm optimization.

\section{Conclusions}
\label{sec:conclusions}

This paper presents a novel surrogate-based optimization framework for WFLO that integrates PCE, a Kriging model, and the EI algorithm. The proposed framework addresses the computational challenges associated with high-fidelity wind farm simulations by significantly reducing the number of function evaluations required for accurate AEP predictions. 
The PCE-based AEP prediction method achieves exceptional accuracy with lower computational expense, significantly lowering the expense of training the ensuing surrogate model.
The Kriging model, combined with a genetic algorithm, is used for surrogate-based optimization, achieving comparable performance to direct optimization at a much reduced computational cost. 
The integration of the EI module enhances the global optimization capability of the framework, allowing it to escape local optima and achieve results that are either nearly identical to or even outperform those obtained through direct optimization. 
The feasibility of the PCE-Kriging framework is demonstrated through four case studies, including the optimization of wind farms with 8, 16, and 32 turbines using low-fidelity wake models, and a high-fidelity case using CFD simulations. 
The results show that the proposed framework is highly effective in optimizing wind farm layouts, significantly reducing computational costs while maintaining or improving the accuracy of AEP predictions.
Overall, the integration of PCE, Kriging, and EI offers a promising approach for future wind farm design, enabling more accurate and efficient optimization of wind farm layouts to maximize energy production and minimize wake losses.

However, the PCE-aided SBO framework remains vulnerable to the curse of dimensionality as the number of turbines increases. 
This issue is a common problem in many machine learning and statistical methods \citep{williams2006gaussian,altman2018curse}.
For the 32-turbine WFLO problem using a low-fidelity wake model studied herein, the wall-clock time for SBO surpasses that of direct optimization, despite the significant reduction in the number of function calls in the former approach.
This inefficiency stems primarily from the computationally expensive process of iteratively updating the Kriging model. 
As the number of turbines grows, the computational time is expected to increase exponentially, further exacerbating the challenge. 
To address this issue, future research should prioritize strategies to mitigate the curse of dimensionality. 
Promising approaches include leveraging dimensionality reduction techniques \citep{bouhlel2016improving} and integrating physics-based knowledge into the optimization process \citep{yang2019physics}, which could enhance both computational efficiency and model accuracy.

\section*{Acknowledgments}

The financial supports from the National Key R\&D Program of China (No. 2023YFE0120000), Guangdong Basic and Applied Basic Research Foundation (No. 2023A1515240054), Program for Intergovernmental International S\&T Cooperation Projects of Shanghai Municipality, China (Nos. 22160710200, 24510711100), National Natural Science Foundation of China (Nos. 42306226) are gratefully acknowledged.

\bibliographystyle{unsrtnat}
\bibliography{reference}

\end{document}